# Intersections of Dragon Curves of the 18th Iteration

Reimund Albers, Zongyi Guo, Huaiyi Guo

Abstract

It is shown that intersections occur in the dragon curve up to an unfolding angle of 95.828°. To this end, two specific edges in the 18th iteration are examined.

## Introduction

The well-known Heighway dragon arises from a paper-folding process in which a strip of paper is folded in half n times, always in the same way. If this strip is then unfolded and the creases are carefully formed into 90° angles, the Heighway dragon of order n is obtained. The distinctive feature of this polygonal path is that no two edges ever lie on top of each other.

If one increases the unfolding angle lager than 90°, the intuitive conjecture was that the polygonal path would then have no intersections at all ([1] M. Mendès France, 1990). In his 2011 book [2], D. Knuth noted that he "noticed in 1969 that 95°-angle folds would lead to paths that cross themselves." In 2006, Albers [3] presented a proof demonstrating the existence of intersecting segments in iteration 10 for angles ranging between 90° and 95.126°.

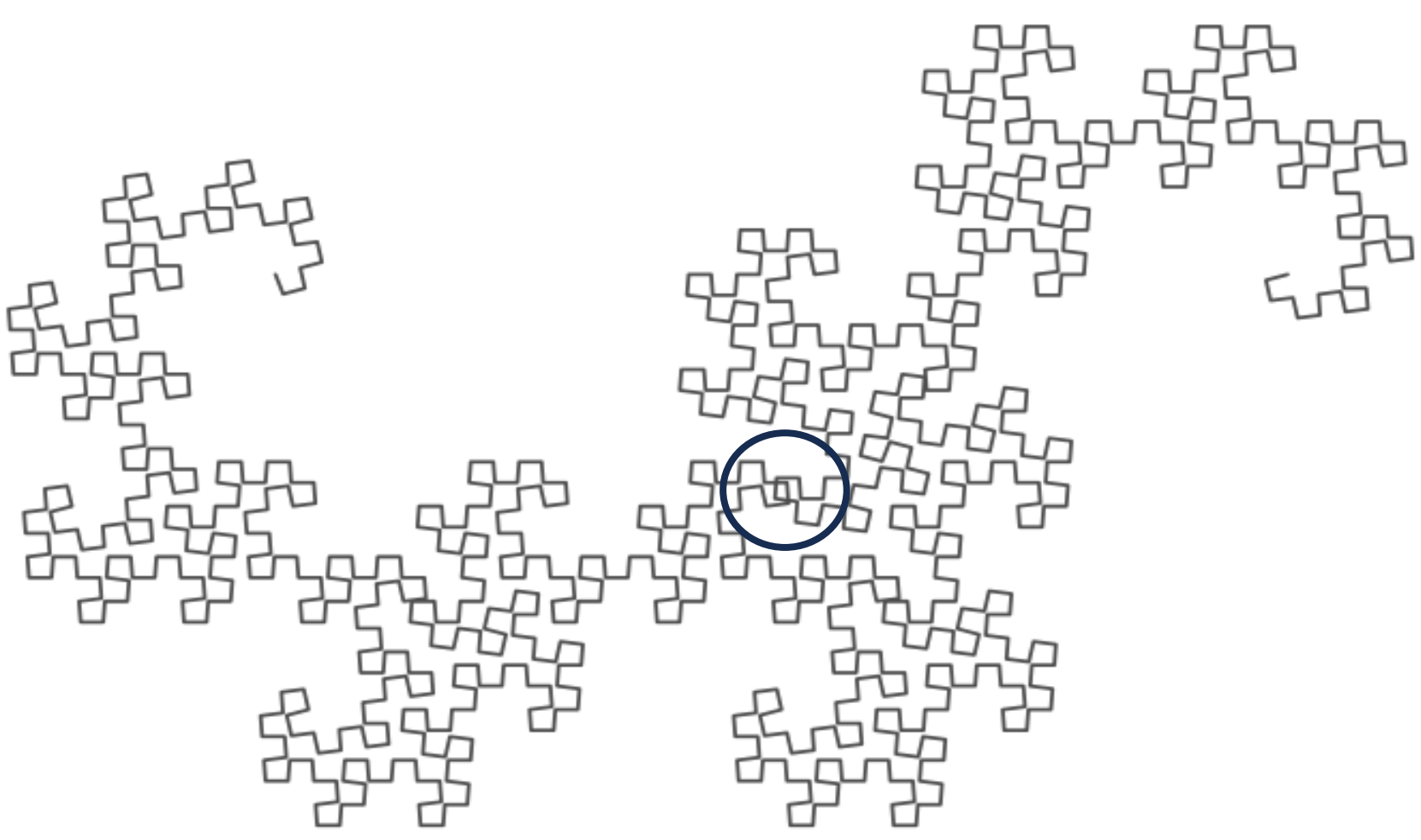

*Figure 1: Intersection in the 10th iteration at 94°*

On the other hand, it has been proven ([4], [5]) that for unfolding angles greater than 98.195° no intersections occur.

This existing body of knowledge is summarized clearly in the following diagram.

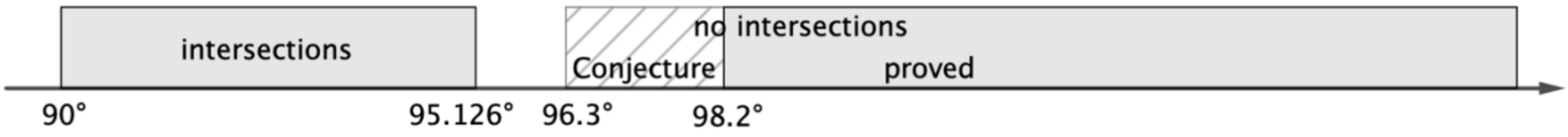

*Figure 2: The previously known bounds for intersections*

In this article, the bound on the unfolding angle up to which intersections can be shown to occur will be raised from 95.126° to 95.828°. To this end, the 18th iteration $Q_{18}$ of the polygonal path is examined. The diagram of the angle bounds then looks as follows:

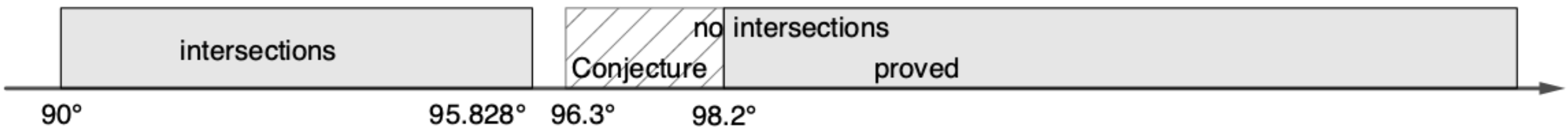

*Figure 3: The new, upper bound for intersections*

## The Paperfolding Polygons $Q_n$ – Construction and Properties

Let us begin with the definition of the objects we will be dealing with.

**Definition 1**

We consider polygons $Q_n$ , $n \in \mathbb{N}_0$ with $2^n + 1$ points $P_{n,k}$ , $k \in \{0, 1, 2, \dots, 2^n\}$ and edges $\overline{P_{n,k}P_{n,k+1}}$ , $k \in \{0, 1, 2, \dots, 2^n - 1\}$. The polygons are called paper-folding polygons or generalized dragon curves of iteration level $n$.

The definition of the $Q_n$ is recursive.

$Q_0$ is the segment $\overline{P_{0,0}P_{0,1}}$ . For the further considerations, we place it in a coordinate system with $P_{0,0} = (0\,, 0)$ und $P_{0,1} = (1\,, 0)$.

$Q_{n+1}$ is constructed from $Q_n$ by

a) Renumbering the vertices $P_{n,k}$ to $P_{n+1,2k}$

b) Replacing each segment $\overline{P_{n,k}P_{n,k+1}}$ by the sides of an isosceles triangle, meeting at $P_{n+1,2k+1}$ under the angle $\alpha$.
   The angles at the base are $\beta$, $\beta = 90° - \frac{\alpha}{2}$.

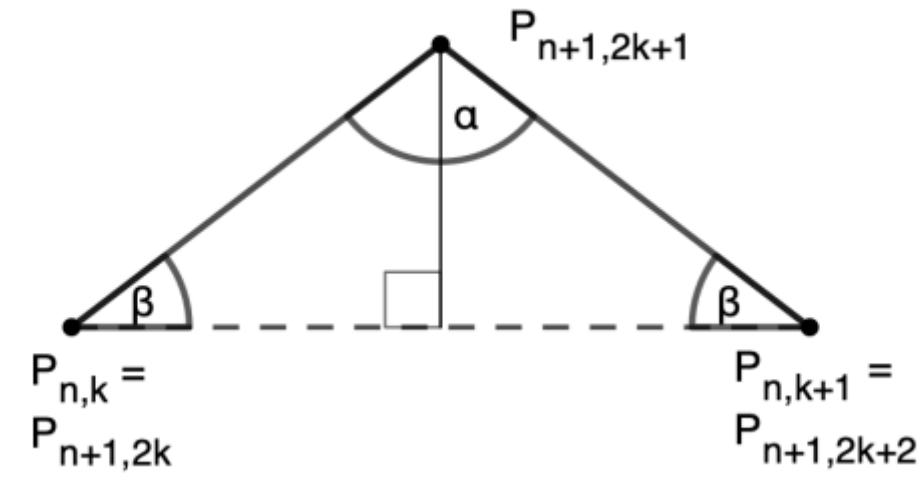

*Figure 4: The angles at the erected triangle*

c) The isosceles triangle is added
   i. to the left side, forming a right turn R at $P_{n+1,2k+1}$, if $\overline{P_{n,k}P_{n,k+1}}$ is from odd *k* to even *k*+1.
   ii. to the right side, forming a left turn L at $P_{n+1,2k+1}$, if $\overline{P_{n,k}P_{n,k+1}}$ is from even *k* to odd *k*+1.

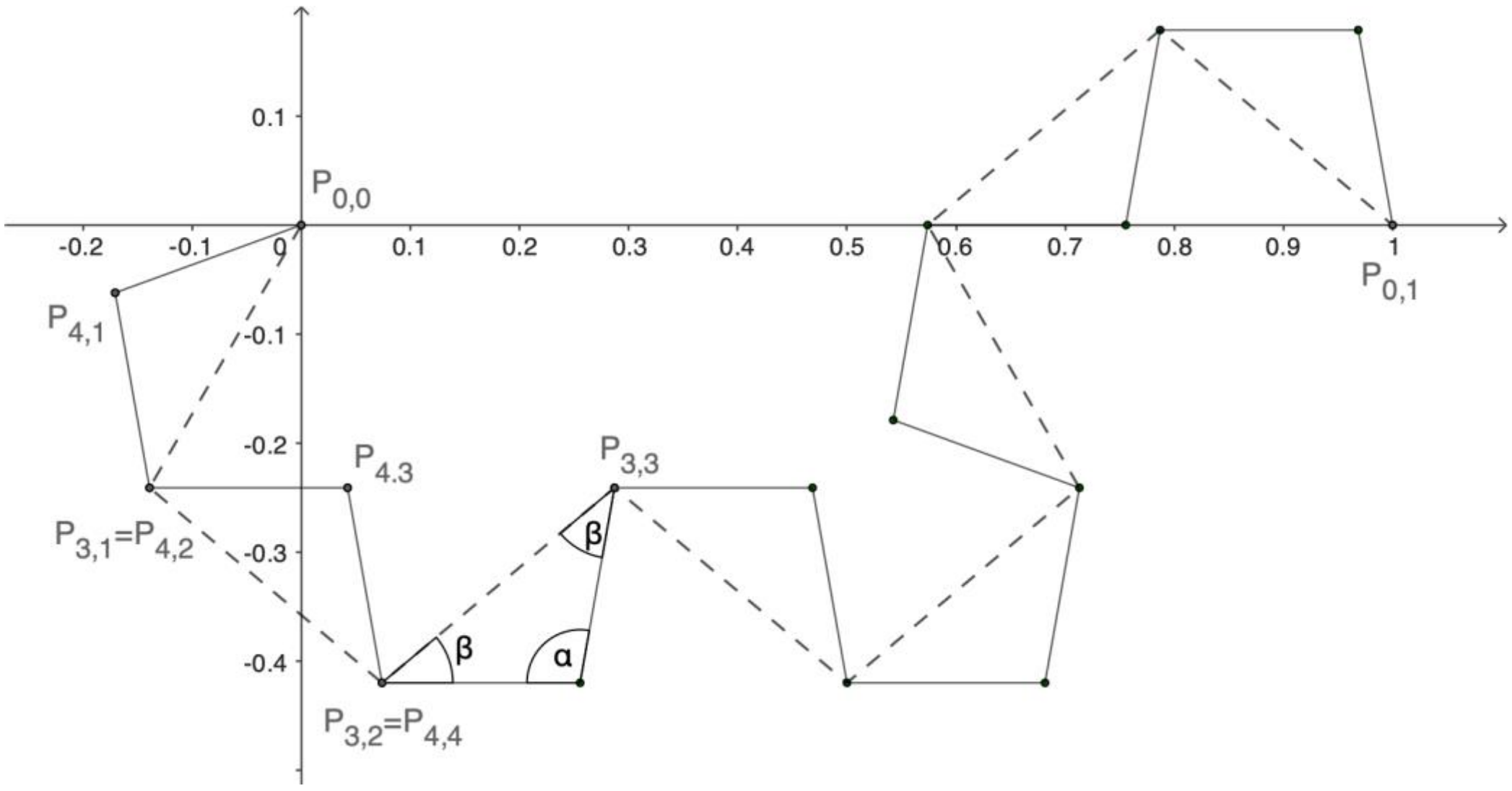

*Figure 5: The construction of $Q_4$ from $Q_3$ (dashed lines)*

**Definition 2**

The function $\sigma$ assigns to each point $P_{n,k}$ , $1 \le k \le 2^n - 1$, its turn direction L or R.

**Lemma 1**

$$\sigma(k) = \begin{cases} \text{L} \ \text{if } k \equiv 1 \ mod\ 4 \\ \text{R} \ \text{if } k \equiv 3 \ mod\ 4 \\ \sigma\left(\frac{k}{2}\right) \ \text{if } k \equiv 0 \ mod\ 2 \end{cases}$$

Proof

The point under consideration is $P_{n,k}$. If $k$ is odd, then $P_{n,k}$ was created as a new vertex on the directed edge $\overrightarrow{P_{n-1,\frac{k-1}{2}}P_{n-1,\frac{k+1}{2}}}$.

a) If $k \equiv 1 \ mod\ 4$, then $k - 1 \equiv 0 \ mod\ 4$, so $\frac{k-1}{2}$ is an even number. Then, according to Definition 1 c) ii., at point $P_{n,k}$ there is a left turn L.

b) If $k \equiv 3 \ mod\ 4$, then $k + 1 \equiv 0 \ mod\ 4$, so $\frac{k+1}{2}$ is an even number and hence $\frac{k-1}{2}$ is an odd number. Then, according to Definition 1 c) i., at point $P_{n,k}$ there is a right turn R.

c) If $k$ is even, then the point $P_{n,k}$ has existed at the prior iteration level, as $P_{n-1,\frac{k}{2}}$.

We consider the construction step from $Q_{n-1}$ to $Q_n$. If $\frac{k}{2}$ is odd, then the edge $\overrightarrow{P_{n-1,\frac{k}{2}-1}P_{n-1,\frac{k}{2}}}$ is an edge from a point with an even number to a point with an odd number. Consequently, the triangle is erected to the right (Definition 1 c) ii.). The edge $\overrightarrow{P_{n-1,\frac{k}{2}}P_{n-1,\frac{k}{2}+1}}$, leading away from the point $P_{n-1,\frac{k}{2}}$, is an edge from a point with an odd number to a point with an even number. Consequently, the triangle is erected to the left (Definition 1 c) i.). (see Figure 6, upper row)

If $\frac{k}{2}$ is even, then the erection of the triangles is exactly reversed (see Figure 6, lower row).

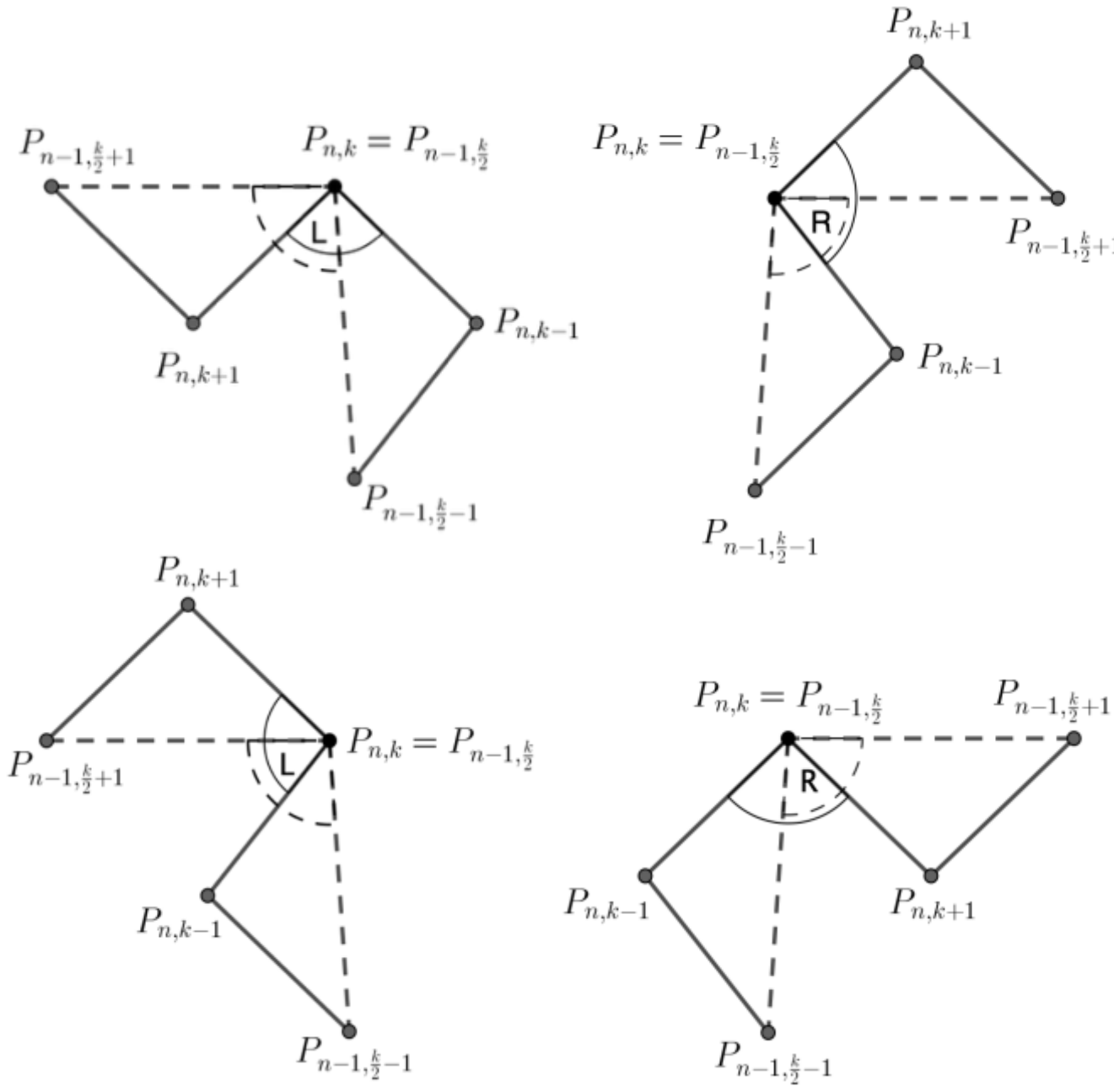

*Figure 6: The orientation of the turns is preserved*

In all four cases, the turn direction is thereby preserved exactly at level $n$ as it already existed at the preceding level $n$-1. □

**Corollary 1**

Lemma 1 yields a very simple algorithm for determining the turn direction at point $P_{n,k}$.
One repeatedly divides the point number $k$ by 2 until the quotient yields an odd number. If its remainder upon division by 4 equals 1, then there is a left turn at point $P_{n,k}$; if the remainder equals 3, then there is a right turn at point $P_{n,k}$

If the point number $k$ is given as a binary number (which is often the case in what follows), the procedure is even simpler and more direct: starting from the right, one looks for the first 1. If the digit before it is a 0, then there is a left turn at point $P_{n,k}$ if it is a 1, then there is a right turn at point $P_{n,k}$.

From the construction rules of Definition 1, the following properties follow immediately:

**Property 1**
a) The edges of the polygon $Q_{n+1}$ are shortened relative to the edges of the polygon $Q_n$ by the factor $0 < q < 1$. $q$ is defined by $q = \frac{1}{2cos\beta}$ (see Fig. 4).
b) Since $|P_{0,0}P_{0,1}| = 1$, the edges of the polygon $Q_n$ have length $q^n$.

**Property 2**
Every point of the polygon $Q_n$ retains its position in the subsequent construction steps. Only its number changes, according to Definition 1 a).
In particular, all polygons start at the point $P_{0,0} = P_{n,0}$ and end at the point $P_{0,1} = P_{n,2^n}$ .

**Definition 3**
Every directed edge $\overrightarrow{P_{n,k}P_{n,k+1}}$ defines an angle to the positive $x$-axis.
We denote this polar angle by $arg\left(\overrightarrow{P_{n,k}P_{n,k+1}}\right)$.

**Property 3**
a) The polar angle of the first edge is $arg\left(\overrightarrow{P_{n,0}P_{n,1}}\right) = -n\beta$.
b) At a left turn L, the polar angle increases by $2\beta$,
i.e. $arg\left(\overrightarrow{P_{n,k+1}P_{n,k+2}}\right) = arg\left(\overrightarrow{P_{n,k}P_{n,k+1}}\right) + 2\beta$
At a right turn R, the polar angle decreases by $2\beta$,
i.e. $arg\left(\overrightarrow{P_{n,k+1}P_{n,k+2}}\right) = arg\left(\overrightarrow{P_{n,k}P_{n,k+1}}\right) - 2\beta$

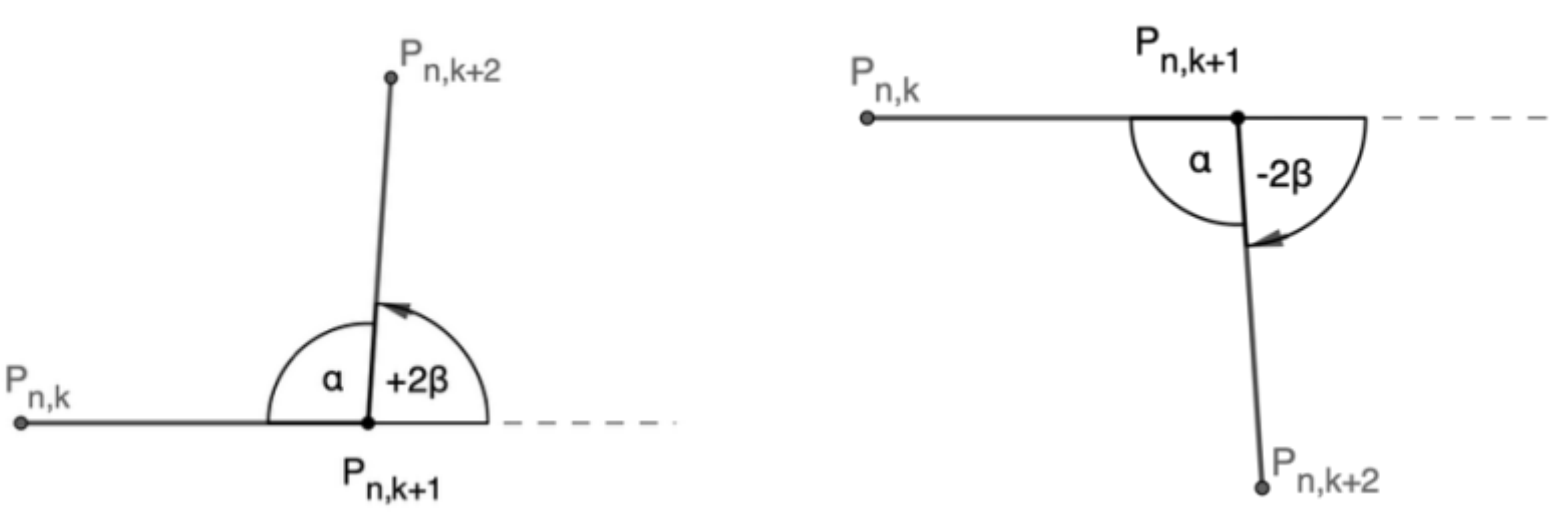

*Figure 7: The change in the polar angle*

**Definition 4**
$run(k)$ determines the number of runs[1] in the binary representation of *k*.
The special definition $run(0) = 0$ holds.
Example: $run(89) = run(1011001_2) = 5$

---
[1] A run in a string is a block of consecutive, identical characters.

**Lemma 2**

$arg(\overrightarrow{P_{n,k}P_{n,k+1}}) = -n\beta + run(k) \cdot 2\beta$ ( see [6] )

Proof

(by mathematical induction on $k$)

Base case:

$k = 0$ Since $run(0) = 0$ we have $arg(\overrightarrow{P_{n,0}P_{n,1}}) = -n\beta$, which is correct according to Property 3 a).

k = 1: By Lemma 1, for all $n$ there is a left turn at the point $P_{n,1}$. By Property 3 b), the argument increases by $2\beta$. so $arg(\overrightarrow{P_{n,1}P_{n,2}}) = arg(\overrightarrow{P_{n,0}P_{n,1}}) + 2\beta = -n\beta + 2\beta$.

With $run(1) = 1$, the formula gives the same result.

Induction hypothesis: For $k \geq 1$, $arg(\overrightarrow{P_{n,k}P_{n,k+1}}) = -n\beta + run(k) \cdot 2\beta$ holds.

Case 1: The binary representation of $k$ is $xx10$ ($xx$ stands for an arbitrary sequence of binary digits)
Then $k + 1 = xx11$. By Corollary 1, there is a right turn at the point $P_{n,k+1}$,

so $arg(\overrightarrow{P_{n,k+1}P_{n,k+2}}) = arg(\overrightarrow{P_{n,k}P_{n,k+1}}) - 2\beta$ and by the induction hypothesis

$$= -n\beta + run(k) \cdot 2\beta - 2\beta = -n\beta + (run(k) - 1) \cdot 2\beta$$

Now here $run(k + 1) = run(k) - 1$, and hence $arg(\overrightarrow{P_{n,k+1}P_{n,k+2}}) = -n\beta + run(k + 1) \cdot 2\beta$.

Case 2: The binary representation of $k$ is $xx1\underbrace{0 \ldots 0}_{a}$ , $a \geq 2$

Then $k + 1 = xx1\underbrace{0 \ldots 0}_{a-1}1$. By Corollary 1, there is a left turn at the point $P_{n,k+1}$, so

$arg(\overrightarrow{P_{n,k+1}P_{n,k+2}}) = arg(\overrightarrow{P_{n,k}P_{n,k+1}}) + 2\beta$ and by the induction hypothesis

$$= -n\beta + run(k) \cdot 2\beta + 2\beta = -n\beta + (run(k) + 1) \cdot 2\beta$$

Now here $run(k + 1) = run(k) + 1$, and hence $arg(\overrightarrow{P_{n,k+1}P_{n,k+2}}) = -n\beta + run(k + 1) \cdot 2\beta$.

Case 3: The binary representation of $k$ is $xx00\underbrace{1 \ldots 1}_{a}$ , $a \geq 1$

Then $k + 1 = xx01\underbrace{0 \ldots 0}_{a}$. By Corollary 1, there is a left turn at the point $P_{n,k+1}$, so

$arg(\overrightarrow{P_{n,k+1}P_{n,k+2}}) = arg(\overrightarrow{P_{n,k}P_{n,k+1}}) + 2\beta$ and by the induction hypothesis

$$= -n\beta + run(k) \cdot 2\beta + 2\beta = -n\beta + (run(k) + 1) \cdot 2\beta$$

Now here $run(k + 1) = run(k) + 1$, and hence $arg(\overrightarrow{P_{n,k+1}P_{n,k+2}}) = -n\beta + run(k + 1) \cdot 2\beta$.

Case 4: The binary representation of $k$ is $xx10\underbrace{1 \ldots 1}_{a}$ , $a \geq 1$

Then $k + 1 = xx11\underbrace{0 \ldots 0}_{a}$. By Corollary 1, there is a right turn at the point $P_{n,k+1}$, so

$arg(\overrightarrow{P_{n,k+1}P_{n,k+2}}) = arg(\overrightarrow{P_{n,k}P_{n,k+1}}) - 2\beta$ and by the induction hypothesis

$$= -n\beta + run(k) \cdot 2\beta - 2\beta = -n\beta + (run(k) - 1) \cdot 2\beta$$

Now here $run(k + 1) = run(k) - 1$, and hence $arg(\overrightarrow{P_{n,k+1}P_{n,k+2}}) = -n\beta + run(k + 1) \cdot 2\beta$.

□

## Preparations for the Calculation

**Algorithm 1** (for the practical calculation of the coordinates of a point $P_{n,k}$)

We calculate a point $P_{n,k}$ by computing the position vector $\overrightarrow{P_{n,0}P_{n,k}}$ as a vector chain that is as short as possible.

$$\overrightarrow{P_{n,0}P_{n,k}} = \overrightarrow{P_{n_1,0}P_{n_1,1}} + \overrightarrow{P_{n_2,k_2}P_{n_2,k_2+1}} + \overrightarrow{P_{n_3,k_3}P_{n_3,k_3+1}} + \cdots + \overrightarrow{P_{n_m,k_m}P_{n_m,k_m+1}} \quad (1)$$

For the iteration levels $n_i$ we have

$$1 \le n_1 < n_2 < n_3 < \cdots < n_m \le n$$

More generally, an edge vector $\overrightarrow{P_{m,k}P_{m,k+1}}$ of a lower iteration level $m < n$ stands for the sum of $2^{n-m}$ edge vectors of level $n$, more precisely

$$\overrightarrow{P_{m,k}P_{m,k+1}} = \sum_{j=0}^{2^{n-m}-1} \overrightarrow{P_{n,k_j}P_{n,k_j+1}} \text{ with } k_j = k \cdot 2^{n-m} + j$$

A justification follows immediately from the construction process and the numbering rule in Definition 1.

For the path from $P_{n,0}$ to the point $P_{n,k}$ the $k$ elementary edge vectors $\overrightarrow{P_{n,i}P_{n,i+1}}$ , $0 \le i \le k-1$ , are combined into powers of two that are as large as possible. This process is thus directly connected to the representation of the number $k$ in the binary system. This gives the following algorithm for calculating the vector chain (1):

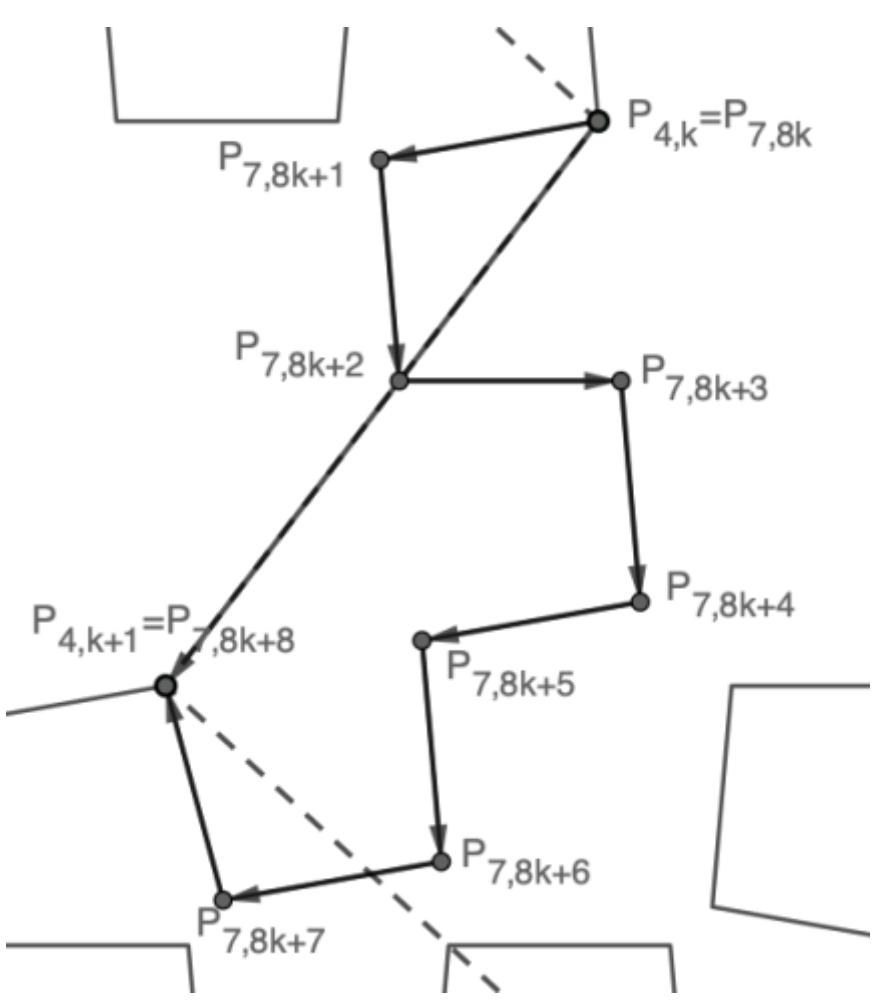

*Figure 8:An edge vector of $Q_4$ stands for the sum of 8 edge vectors of $Q_7$*

Step 1: The number $k$ is represented in the binary system, $k = \sum_{i=0}^{n} b_i 2^i$ , $b_i \in \{0,1\}$.

Step 2: One determines, from left to right (i.e., starting with the highest power of two), the indices $n-1 \ge i_1 > i_2 > i_3 > \cdots > i_m \ge 0$, for which $b_{i_j} = 1$.

Step 3: The iteration levels to be used for the vector chain are then $n_j = n - i_j$ , $1 \le j \le m$. Hence $1 \le n_1 < n_2 < n < \cdots < n_m \le n$.

Step 4: The point numbers $k_j$ belonging to the individual iteration levels $n_j$ are obtained recursively by: $k_1 = 0$ und $k_{j+1} = (k_j + 1) \cdot 2^{n_{j+1}-n_j}$

**Example:** Given the point $P_{8,89}$, i.e. $n = 8, k = 89$.

Step 1: $k = 89 = 1011001_2$

Step 2: The digit 1 occurs at the positions $i_1 = 6$ , $i_2 = 4$ , $i_3 = 3$ , $i_4 = 0$ $(m = 4)$

Step 3: The iteration levels to be used are then

$$n_1 = n - i_1 = 2 \text{ , } n_2 = n - i_2 = 4 \text{ , } n_3 = n - i_3 = 5 \text{ , } n_4 = n - i_4 = 8$$

Step 4: The point numbers for the starting points of the individual sub-vectors are

$$k_1 = 0$$
$$k_2 = (k_1 + 1) \cdot 2^{n_2-n_1} = 1 \cdot 2^{4-2} = 4$$
$$k_3 = (k_2 + 1) \cdot 2^{n_3-n_2} = 5 \cdot 2^{5-4} = 10$$
$$k_4 = (k_3 + 1) \cdot 2^{n_4-n_3} = 11 \cdot 2^{8-5} = 88$$

Thus, the vector chain is: $\overrightarrow{P_{8,0}P_{8,89}} = \overrightarrow{P_{2,0}P_{2,1}} + \overrightarrow{P_{4,4}P_{4,5}} + \overrightarrow{P_{5,10}P_{5,11}} + \overrightarrow{P_{8,88}P_{8,89}}$

**Algorithm 1'**: (further shortening of the vector chain)
If, in the binary representation of $k$, one obtains a longer (at least 3) run of consecutive ones, the number of positions at which the binary representation contains no 0 can be reduced.
For this, one uses the identity $2^n + 2^{n-1} + 2^{n-2} + \cdots + 2^m = 2^{n+1} - 2^m$.
A binary representation ... 0 $\underbrace{111 \ldots 11}_{a\ digits}$ 0 ... then becomes ... 1 $\underbrace{000 \ldots 0\overline{1}}_{a\ digits}$ 0 ..., where $\overline{1}$ stands for -1.
The algorithm, adapted to this more general form, is then:
Step 1: a) The number $k$ is represented in the binary system, $k = \sum_{i=0}^{n} b_i 2^i \ , b_i \in \{0,1\}$.
b) Longer runs of ones are converted using $\overline{1}$.
Step 2: One determines, from the right (starting with the highest power of two), the indices $n-1 \geq i_1 > i_2 > i_3 > \cdots > i_m \geq 0$, for which $\left|b_{i_j}\right| = 1$.
Step 3: The iteration levels to be used for the vector chain are then $n_j = n - i_j \ , 1 \leq j \leq m$. Hence $1 \leq n_1 < n_2 < n < \cdots < n_m \leq n$.
Step 4: The point numbers $k_j$ belonging to the individual iteration levels $n_j$ are obtained recursively by: $k_1 = 0$ und $k_{j+1} = \left(k_j + b_{i_j}\right) \cdot 2^{n_{j+1}-n_j}$
In addition, the vector chain must be modified. From the starting point $P_{n_j,k_j}$ of a vector, one moves forward or backward, depending on whether the corresponding digit is +1 or -1.

$$\overrightarrow{P_{n,0}P_{n,k}} = \overrightarrow{P_{n_1,0}P_{n_1,1}} + \overrightarrow{P_{n_2,k_2}P_{n_2,k_2+b_{\iota_2}}} + \overrightarrow{P_{n_3,k_3}P_{n_3,k_3+b_{\iota_3}}} + \cdots + \overrightarrow{P_{n_m,k_m}P_{n_m,k_m+b_{\iota_m}}} \quad (2)$$

**Example:** Given the point $P_{8,189}$, i.e. $n = 8, k = 189$.
Step 1: a) $k = 189 = 10111101_2$
b) $k = 189 = 1100\overline{1}101_2$
Step 2: A digit $\pm 1$ occurs at the positions $i_1 = 7 \ , i_2 = 6 \ , i_3 = 2 \ , i_4 = 0 \ \ (m = 4)$
and the binary digits are $b_{i_1} = 1 \ , b_{i_2} = 1 \ , b_{i_3} = -1 \ , b_{i_4} = 1$ ,
Step 3: The iteration levels to be used are then
$n_1 = n - i_1 = 1 \ , n_2 = n - i_2 = 2 \ , n_3 = n - i_3 = 6 \ , n_4 = n - i_4 = 8$
Step 4: The point numbers for the starting points of the individual sub-vectors are
$k_1 = 0$
$k_2 = (k_1 + 1) \cdot 2^{n_2 - n_1} = 1 \cdot 2^{2-1} = 2$
$k_3 = (k_2 + 1) \cdot 2^{n_3 - n_2} = 3 \cdot 2^{6-2} = 48$
$k_4 = (k_3 - 1) \cdot 2^{n_4 - n_3} = 47 \cdot 2^{8-6} = 188$
Thus, the vector chain is: $\overrightarrow{P_{8,0}P_{8,189}} = \overrightarrow{P_{1,0}P_{1,1}} + \overrightarrow{P_{2,2}P_{2,3}} + \overrightarrow{P_{6,48}P_{6,47}} + \overrightarrow{P_{8,188}P_{8,189}}$

Since in the vector chain the vectors $\overrightarrow{P_{n,k}P_{n,k+1}}$ are added, these must be given in Cartesian coordinates. The calculation from polar coordinates is carried out in the usual way

$$\overrightarrow{P_{n,k}P_{n,k+1}} = \left|P_{n,k}P_{n,k+1}\right| \begin{pmatrix} \cos(\arg(\overrightarrow{P_{n,k}P_{n,k+1}})) \\ \sin(\arg(\overrightarrow{P_{n,k}P_{n,k+1}})) \end{pmatrix}$$

The coordinates calculated in this way are to be expressed in terms of the variable $q$. For this reason, a compilation follows for the conversion of some trigonometric functions of the angle $\beta$.

Compilation of Formulas for Conversion

$$q = \frac{1}{2\cos\beta} \Rightarrow \cos\beta = \frac{1}{2q}\ , \qquad \sin\beta = \sqrt{1-\cos^2\beta} = \frac{1}{2q}\sqrt{4q^2-1}\ \text{, da } 0 \le \beta < 90°$$

$$\cos 2\beta = 2\cos^2\beta - 1 = \frac{1}{2q^2}\,(1-2q^2)\,, \qquad \sin 2\beta = 2\sin\beta\cos\beta = \frac{1}{2q^2}\sqrt{4q^2-1}$$

$$\cos 3\beta = 4\cos^3\beta - 3\cos\beta = \frac{1}{2q^3}\,(1-3q^2)$$

$$\sin 3\beta = \sin\beta\,(4\cos^2\beta - 1) = \frac{1}{2q^3}\sqrt{4q^2-1}\,(1-q^2)$$

$$\cos 4\beta = 2\cos^2(2\beta) - 1 = \frac{1}{2q^4}\,(1-4q^2+2q^4)$$

$$\sin 4\beta = 2\sin 2\beta\cos 2\beta = \frac{1}{2q^4}\sqrt{4q^2-1}\,(1-2q^2)$$

$$\cos 5\beta = 16\cos^5\beta - 20\cos^3\beta + 5\cos\beta = \frac{1}{2q^5}\,(1-5q^2+5q^4)$$

$$\sin 5\beta = \sin\beta\,(16\cos^4\beta - 12\cos^2\beta + 1) = \frac{1}{2q^5}\sqrt{4q^2-1}\,(1-3q^2+q^4)$$

## Calculating the Intersection Point

**Theorem 1**
For an unfolding angle from 95,752° to 95,828°, intersections occur in the paper-folding polygon $Q_{18}$ between the edges $\overrightarrow{P_{18,114167}P_{18,114168}}$ and $\overrightarrow{P_{18,179694}P_{18,179695}}$.

Proof
a) Calculation of the coordinates of $P_{18,114168}$
Following Algorithm 1, we first determine the binary representation of 114168.

| Exponent for 2 | 17 | 16 | 15 | 14 | 13 | 12 | 11 | 10 | 9 | 8 | 7 | 6 | 5 | 4 | 3 | 2 | 1 | 0 |
|---|---|---|---|---|---|---|---|---|---|---|---|---|---|---|---|---|---|---|
| Binary digits | 0 | 1 | 1 | 0 | 1 | 1 | 1 | 1 | 0 | 1 | 1 | 1 | 1 | 1 | 1 | 0 | 0 | 0 |
| Digits using -1 | 0 | 1 | 1 | 1 | 0 | 0 | 0 | -1 | 1 | 0 | 0 | 0 | 0 | 0 | -1 | 0 | 0 | 0 |
| Iteration level | 1 | 2 | 3 | 4 | 5 | 6 | 7 | 8 | 9 | 10 | 11 | 12 | 13 | 14 | 15 | 16 | 17 | 18 |

Twice a longer run of ones occurs, so it is worthwhile to use the simplification via the modified Algorithm 1' (using −1) here.
This gives the following vector sum for the position vector $\overrightarrow{P_{18,0}P_{18,114168}}$:

$$\overrightarrow{P_{18,0}P_{18,114168}} = \overrightarrow{P_{2,0}P_{2,1}} + \overrightarrow{P_{3,2}P_{3,3}} + \overrightarrow{P_{4,6}P_{4,7}} + \overrightarrow{P_{8,112}P_{8,111}} + \overrightarrow{P_{9,222}P_{9,223}} + \overrightarrow{P_{15,14272}P_{15,14271}}$$

The endpoint $P_{15,14271}$ of the vector chain coincides with the target point $P_{18,114168}$.

Calculation of the individual sub-vectors

$\overrightarrow{P_{2,0}P_{2,1}}$
$\left|\overrightarrow{P_{2,0}P_{2,1}}\right| = q^2$ (Property 1 b) )

$arg\left(\overrightarrow{P_{2,0}P_{2,1}}\right) = -2\beta$ (Property 3 a) )

$$\overrightarrow{P_{2,0}P_{2,1}} = q^2 \begin{pmatrix} cos(-2\beta) \\ sin(-2\beta) \end{pmatrix} = \begin{pmatrix} q^2\ cos(2\beta) \\ -q^2\ sin(2\beta) \end{pmatrix} = \begin{pmatrix} \frac{1}{2}\,(1-2q^2) \\ -\frac{1}{2}\sqrt{4q^2-1} \end{pmatrix}$$

$\overrightarrow{P_{3,2}P_{3,3}}$

$\left|\overrightarrow{P_{3,2}P_{3,3}}\right| = q^3$

$k = 2 = 10_2$ , hence $run(2) = 2$

$arg\left(\overrightarrow{P_{3,2}P_{3,3}}\right) = -3\beta + 2 \cdot 2\beta = \beta$

$$\overrightarrow{P_{3,2}P_{3,3}} = q^3 \begin{pmatrix} cos(\beta) \\ sin(\beta) \end{pmatrix} = \begin{pmatrix} \frac{1}{2}q^2 \\ \frac{1}{2}q^2\sqrt{4q^2-1} \end{pmatrix}$$

$\overrightarrow{P_{4,6}P_{4,7}}$

$\left|\overrightarrow{P_{4,6}P_{4,7}}\right| = q^4$

$k = 6 = 110_2$, hence $run(6) = 2$

$arg\left(\overrightarrow{P_{4,6}P_{4,7}}\right) = -4\beta + 2 \cdot 2\beta = 0°$

$$\overrightarrow{P_{4,6}P_{4,7}} = q^4 \begin{pmatrix} cos(0°) \\ sin(0°) \end{pmatrix} = \begin{pmatrix} q^4 \\ 0 \end{pmatrix}$$

$\overrightarrow{P_{8,112}P_{8,111}} = -\overrightarrow{P_{8,111}P_{8,112}}$

$\left|\overrightarrow{P_{8,111}P_{8,112}}\right| = q^8$

$k = 111 = 1101111_2$ , hence $run(111) = 3$

$arg\left(\overrightarrow{P_{8,111}P_{8,112}}\right) = -8\beta + 3 \cdot 2\beta = -2\beta$

$$\overrightarrow{P_{8,112}P_{8,111}} = -\overrightarrow{P_{8,111}P_{8,112}} = -q^8 \begin{pmatrix} cos(-2\beta) \\ sin(-2\beta) \end{pmatrix} = \begin{pmatrix} -\frac{1}{2}\ q^6(1-2q^2) \\ \frac{1}{2}q^6\sqrt{4q^2-1} \end{pmatrix}$$

$\overrightarrow{P_{9,222}P_{9,223}}$

$\left|\overrightarrow{P_{9,222}P_{9,223}}\right| = q^9$

$k = 222 = 11011110_2$ , hence $run(222) = 4$

$arg\left(\overrightarrow{P_{9,222}P_{9,223}}\right) = -9\beta + 4 \cdot 2\beta = -\beta$

$$\overrightarrow{P_{9,222}P_{9,223}} = q^9 \begin{pmatrix} cos(-\beta) \\ sin(-\beta) \end{pmatrix} = \begin{pmatrix} \frac{1}{2}q^8 \\ -\frac{1}{2}q^8\sqrt{4q^2-1} \end{pmatrix}$$

$\overrightarrow{P_{15,14272}P_{15,14271}} = -\overrightarrow{P_{15,14271}P_{15,14272}}$

$\left|\overrightarrow{P_{15,14271}P_{15,14272}}\right| = q^{15}$

$k = 14271 = 11011110111111_2$ , hence $run(14271) = 5$

$arg\left(\overrightarrow{P_{15,14271}P_{15,14272}}\right) = -15\beta + 5 \cdot 2\beta = -5\beta$

$$\overrightarrow{P_{15,14272}P_{15,14271}} = -\overrightarrow{P_{15,14271}P_{15,14272}} = -q^{15}\begin{pmatrix} cos(-5\beta) \\ sin(-5\beta) \end{pmatrix} = \begin{pmatrix} -\frac{1}{2}\,q^{10}(1-5q^2+5q^4) \\ \frac{1}{2}q^{10}\sqrt{4q^2-1}\,(1-3q^2+q^4) \end{pmatrix}$$

Adding all the sub-vectors gives

$$\overrightarrow{P_{18,0}P_{18,114168}} = \frac{1}{2}\begin{pmatrix} 1-q^2+2q^4-q^6+3q^8-q^{10}+5q^{12}-5q^{14} \\ (-1+q^2+q^6-q^8+q^{10}-3q^{12}+q^{14})\sqrt{4q^2-1} \end{pmatrix}$$

b) Calculation of the vector $\overrightarrow{P_{18,114167}P_{18,114168}}$

$|\overrightarrow{P_{18,114167}P_{18,114168}}| = q^{18}$

$k = 114167 = 11011110111110111_2$ folglich $run(114167) = 7$

$arg(\overrightarrow{P_{18,114167}P_{18,114168}}) = -18\beta + 7\cdot 2\beta = -4\beta$

$$\overrightarrow{P_{18,114167}P_{18,114168}} = q^{18}\begin{pmatrix} cos(-4\beta) \\ sin(-4\beta) \end{pmatrix} = \begin{pmatrix} \frac{1}{2}\,q^{14}(1-4q^2+2q^4) \\ -\frac{1}{2}q^{14}\sqrt{4q^2-1}\,(1-2q^2) \end{pmatrix}$$

The position vector for all points $P$ on the segment $\overline{P_{18,114167}P_{18,114168}}$ is obtained by subtracting a multiple of the vector $\overrightarrow{P_{18,114167}P_{18,114168}}$ from the vector $\overrightarrow{P_{18,0}P_{18,114168}}$.

$$\overrightarrow{P_{18,0}P} = \overrightarrow{P_{18,0}P_{18,114168}} - s\,\overrightarrow{P_{18,114167}P_{18,114168}}\,,0 \le s \le 1$$

$$\overrightarrow{P_{18,0}P} = \frac{1}{2}\begin{pmatrix} 1-q^2+2q^4-q^6+3q^8-q^{10}+5q^{12}-5q^{14} \\ (-1+q^2+q^6-q^8+q^{10}-3q^{12}+q^{14})\sqrt{4q^2-1} \end{pmatrix} - \frac{1}{2}\,s\,q^{14}\begin{pmatrix} (1-4q^2+2q^4) \\ -(1-2q^2)\sqrt{4q^2-1} \end{pmatrix}$$

c) Calculation of the coordinates of $P_{18,179694}$

Following Algorithm 1, we first determine the binary representation of 179694.

| Exponent of 2 | 17 | 16 | 15 | 14 | 13 | 12 | 11 | 10 | 9 | 8 | 7 | 6 | 5 | 4 | 3 | 2 | 1 | 0 |
|---|---|---|---|---|---|---|---|---|---|---|---|---|---|---|---|---|---|---|
| Binary digits | 1 | 0 | 1 | 0 | 1 | 1 | 1 | 1 | 0 | 1 | 1 | 1 | 1 | 0 | 1 | 1 | 1 | 0 |
| Digits using -1 | 1 | 0 | 1 | 1 | 0 | 0 | 0 | -1 | 1 | 0 | 0 | 0 | -1 | 1 | 0 | 0 | -1 | 0 |
| Iteration level | 1 | 2 | 3 | 4 | 5 | 6 | 7 | 8 | 9 | 10 | 11 | 12 | 13 | 14 | 15 | 16 | 17 | 18 |

Three times a longer run of ones occurs, so it is worthwhile to use the simplification via the modified Algorithm 1' (using −1) here.

This gives the following vector sum for the position vector $\overrightarrow{P_{18,0}P_{18,179694}}$:

$$\overrightarrow{P_{18,0}P_{18,179694}} = \overrightarrow{P_{1,0}P_{1,1}} + \overrightarrow{P_{3,4}P_{3,5}} + \overrightarrow{P_{4,10}P_{4,11}} + \overrightarrow{P_{8,176}P_{8,175}} + \overrightarrow{P_{9,350}P_{9,351}} + \overrightarrow{P_{13,5616}P_{13,5615}} + \overrightarrow{P_{14,11230}P_{14,11231}} + \overrightarrow{P_{17,89848}P_{17,89847}}$$

The endpoint $P_{17,89847}$ of the vector chain coincides with the target point $P_{18,179694}$.

Calculation of the individual sub-vectors

$\overrightarrow{P_{1,0}P_{1,1}}$

$|\overrightarrow{P_{1,0}P_{1,1}}| = q$ (Property 1 b) )

$arg(\overrightarrow{P_{1,0}P_{1,1}}) = -\beta$ (Property 3 a) )

$$\overrightarrow{P_{1,0}P_{1,1}} = q\begin{pmatrix} cos(-\beta) \\ sin(-\beta) \end{pmatrix} = \begin{pmatrix} q\ cos(\beta) \\ -q\ sin(\beta) \end{pmatrix} = \begin{pmatrix} \frac{1}{2} \\ -\frac{1}{2}\sqrt{4q^2-1} \end{pmatrix}$$

$\overrightarrow{P_{3,4}P_{3,5}}$
$\left|\overrightarrow{P_{3,4}P_{3,5}}\right| = q^3$
$k = 4 = 100_2$ , hence $run(2) = 2$
$arg\left(\overrightarrow{P_{3,4}P_{3,5}}\right) = -3\beta + 2 \cdot 2\beta = \beta$

$$\overrightarrow{P_{3,4}P_{3,5}} = q^3\begin{pmatrix} cos(\beta) \\ sin(\beta) \end{pmatrix} = \begin{pmatrix} \frac{1}{2}q^2 \\ \frac{1}{2}q^2\sqrt{4q^2-1} \end{pmatrix}$$

$\overrightarrow{P_{4,10}P_{4,11}}$
$\left|\overrightarrow{P_{4,10}P_{4,11}}\right| = q^4$
$k = 10 = 1010_2$ , hence $run(10) = 4$
$arg\left(\overrightarrow{P_{4,10}P_{4,11}}\right) = -4\beta + 4 \cdot 2\beta = 4\beta$

$$\overrightarrow{P_{4,10}P_{4,11}} = q^4\begin{pmatrix} cos(4\beta) \\ sin(4\beta) \end{pmatrix} = \begin{pmatrix} \frac{1}{2}(1-4q^2+2q^4) \\ \frac{1}{2}\sqrt{4q^2-1}\ (1-2q^2) \end{pmatrix}$$

$\overrightarrow{P_{8,176}P_{8,175}} = -\overrightarrow{P_{8,175}P_{8,176}}$
$\left|\overrightarrow{P_{8,175}P_{8,176}}\right| = q^8$
$k = 175 = 10101111_2$ , hence $run(175) = 5$
$arg\left(\overrightarrow{P_{8,175}P_{8,176}}\right) = -8\beta + 5 \cdot 2\beta = 2\beta$

$$\overrightarrow{P_{8,176}P_{8,175}} = -\overrightarrow{P_{8,175}P_{8,176}} = -q^8\begin{pmatrix} cos(2\beta) \\ sin(2\beta) \end{pmatrix} = \begin{pmatrix} -\frac{1}{2}\ q^6(1-2q^2) \\ -\frac{1}{2}q^6\sqrt{4q^2-1} \end{pmatrix}$$

$\overrightarrow{P_{9,350}P_{9,351}}$
$\left|\overrightarrow{P_{9,350}P_{9,351}}\right| = q^9$
$k = 350 = 101011110_2$ , hence $run(2) = 6$
$arg\left(\overrightarrow{P_{9,350}P_{9,351}}\right) = -9\beta + 6 \cdot 2\beta = 3\beta$

$$\overrightarrow{P_{9,350}P_{9,351}} = q^9\begin{pmatrix} cos(3\beta) \\ sin(3\beta) \end{pmatrix} = \begin{pmatrix} \frac{1}{2}q^6(1-3q^2) \\ \frac{1}{2}q^6\sqrt{4q^2-1}\ (1-q^2) \end{pmatrix}$$

$\overrightarrow{P_{13,5616}P_{13,5615}} = -\overrightarrow{P_{13,5615}P_{13,5616}}$
$\left|\overrightarrow{P_{13,5615}P_{13,5616}}\right| = q^{13}$
$k = 5615 = 1010111101111_2$ , hence $run(5615) = 7$
$arg\left(\overrightarrow{P_{13,5615}P_{13,5616}}\right) = -13\beta + 7 \cdot 2\beta = \beta$

$$\overrightarrow{P_{13,5616}P_{13,5615}} = -\overrightarrow{P_{13,5615}P_{13,5616}} = -q^{13}\begin{pmatrix} cos(\beta) \\ sin(\beta) \end{pmatrix} = \begin{pmatrix} -\frac{1}{2}\,q^{12} \\ -\frac{1}{2}q^{12}\sqrt{4q^2-1} \end{pmatrix}$$

$\overrightarrow{P_{14,11230}P_{14,11231}}$
$\left|\overrightarrow{P_{14,11230}P_{14,11231}}\right| = q^{14}$
$k = 11230 = 10101111011110_2$ , hence $run(11230) = 8$
$arg\left(\overrightarrow{P_{14,11230}P_{14,11231}}\right) = -14\beta + 8 \cdot 2\beta = 2\beta$

$$\overrightarrow{P_{14,11230}P_{14,11231}} = q^{14}\begin{pmatrix} cos(2\beta) \\ sin(2\beta) \end{pmatrix} = \begin{pmatrix} \frac{1}{2}q^{12}(1-2q^2) \\ \frac{1}{2}q^{12}\sqrt{4q^2-1} \end{pmatrix}$$

$\overrightarrow{P_{17,89848}P_{17,89847}} = -\overrightarrow{P_{17,89847}P_{17,89848}}$
$\left|\overrightarrow{P_{17,89847}P_{17,89848}}\right| = q^{17}$
$k = 89847 = 10101111011110111_2$ , hence $run(89847) = 9$
$arg\left(\overrightarrow{P_{17,89847}P_{17,89848}}\right) = -17\beta + 9 \cdot 2\beta = \beta$

$$\overrightarrow{P_{17,89848}P_{17,89847}} = -\overrightarrow{P_{17,89847}P_{17,89848}} = -q^{17}\begin{pmatrix} cos(\beta) \\ sin(\beta) \end{pmatrix} = \begin{pmatrix} -\frac{1}{2}\,q^{16} \\ -\frac{1}{2}q^{16}\sqrt{4q^2-1} \end{pmatrix}$$

Adding all the sub-vectors gives

$$\overrightarrow{P_{18,0}P_{18,179694}} = \frac{1}{2}\begin{pmatrix} 2-3q^2+2q^4-q^8-2q^{14}-q^{16} \\ (-q^2-q^8-q^{16})\sqrt{4q^2-1} \end{pmatrix}$$

d) Calculation of the vector $\overrightarrow{P_{18,179694}P_{18,179695}}$
$\left|\overrightarrow{P_{18,179694}P_{18,179695}}\right| = q^{18}$
$k = 179694 = 101011110111101110_2$ , hence $run(179694) = 10$
$arg\left(\overrightarrow{P_{18,179694}P_{18,179695}}\right) = -18\beta + 10 \cdot 2\beta = 2\beta$

$$\overrightarrow{P_{18,179694}P_{18,179695}} = q^{18}\begin{pmatrix} cos(2\beta) \\ sin(2\beta) \end{pmatrix} = \begin{pmatrix} \frac{1}{2}\,q^{16}(1-2q^2) \\ \frac{1}{2}q^{16}\sqrt{4q^2-1} \end{pmatrix}$$

The position vector for all points Q of the segment $\overline{P_{18,179694}P_{18,179695}}$ is obtained by adding a multiple of the vector $\overrightarrow{P_{18,179694}P_{18,179695}}$ to the vector $\overrightarrow{P_{18,0}P_{18,179694}}$.

$$\overrightarrow{P_{18,0}Q} = \overrightarrow{P_{18,0}P_{18,179694}} + t\,\overrightarrow{P_{18,179694}P_{18,179695}},\ 0 \le t \le 1$$

$$\overrightarrow{P_{18,0}Q} = \frac{1}{2}\begin{pmatrix} 2-3q^2+2q^4-q^8-2q^{14}-q^{16} \\ (-q^2-q^8-q^{16})\sqrt{4q^2-1} \end{pmatrix} + \frac{1}{2}t\,q^{16}\begin{pmatrix} 1-2q^2 \\ \sqrt{4q^2-1} \end{pmatrix}$$

e) Calculation of the intersection point of the two lines $P_{18,114167}P_{18,114168}$ and $P_{18,179694}P_{18,179695}$
The condition for the two edges $\overline{P_{18,114167}P_{18,114168}}$ and $\overline{P_{18,179694}P_{18,179695}}$ to have a common intersection point is

$$\overrightarrow{P_{18,0}P} = \overrightarrow{P_{18,0}Q} \text{ and } 0 \le s,t \le 1.$$

This gives the following system (2x2) of linear equations:

$$\frac{1}{2}(1-q^2+2q^4-q^6+3q^8-q^{10}+5q^{12}-5q^{14})-\frac{1}{2}\,s\,q^{14}(1-4q^2+2q^4)$$
$$=\frac{1}{2}\left(2-3q^2+2q^4-q^8-2q^{14}-q^{16}\right)+\frac{1}{2}t\,q^{16}(1-2q^2)$$
$$\frac{1}{2}(-1+q^2+q^6-q^8+q^{10}-3q^{12}+q^{14})\sqrt{4q^2-1}+\frac{1}{2}\,s\,q^{14}(1-2q^2)\sqrt{4q^2-1}$$
$$=\frac{1}{2}(-q^2-q^8-q^{16})\sqrt{4q^2-1}+\frac{1}{2}t\,q^{16}\sqrt{4q^2-1}$$

The solutions for the variables $s$ and $t$ are:

$$s=\frac{-1+2q^2-q^4+3q^6-q^8+5q^{10}-5q^{12}+q^{14}+q^{16}}{q^{12}(1-4q^2+3q^4)}$$
$$t=\frac{-1+5q^2-7q^4+2q^6+q^8-3q^{10}+q^{14}-q^{16}-2q^{18}+q^{20}}{q^{16}(1-4q^2+3q^4)}$$

The range of $q$ for which the conditions $0\le s,t\le 1$ are satisfied is determined numerically.

| alpha | beta | q | s | t |
|---|---|---|---|---|
| 95,750 | 42,1250 | 0,674142 | -0,0228 | 0,8639 |
| 95,751 | 42,1245 | 0,674137 | -0,0098 | 0,8564 |
| 95,752 | 42,1240 | 0,674132 | 0,0031 | 0,8490 |
| 95,753 | 42,1235 | 0,674126 | 0,0161 | 0,8415 |
| 95,754 | 42,1230 | 0,674121 | 0,0290 | 0,8340 |
| 95,755 | 42,1225 | 0,674116 | 0,0420 | 0,8266 |
| 95,756 | 42,1220 | 0,674110 | 0,0549 | 0,8191 |
| 95,757 | 42,1215 | 0,674105 | 0,0679 | 0,8117 |
| ... | ... | ... | ... | ... |
| 95,822 | 42,0890 | 0,673760 | 0,9170 | 0,3252 |
| 95,823 | 42,0885 | 0,673754 | 0,9302 | 0,3177 |
| 95,824 | 42,0880 | 0,673749 | 0,9434 | 0,3102 |
| 95,825 | 42,0875 | 0,673744 | 0,9565 | 0,3027 |
| 95,826 | 42,0870 | 0,673738 | 0,9697 | 0,2952 |
| 95,827 | 42,0865 | 0,673733 | 0,9829 | 0,2877 |
| 95,828 | 42,0860 | 0,673728 | 0,9961 | 0,2802 |
| 95,829 | 42,0855 | 0,673722 | 1,0092 | 0,2727 |
| 95,830 | 42,0850 | 0,673717 | 1,0224 | 0,2651 |

*Table 1*

For $95{,}752° \le \alpha \le 95{,}828°$, the intersection point of $\overline{P_{18,114167}P_{18,114168}}$ and $\overline{P_{18,179694}P_{18,179695}}$ is an interior point of both segments, so that $Q_{18}$ has an intersection at least at this location. □